\documentclass{article}

\usepackage{amsmath}
\usepackage{url}
\usepackage[T2A]{fontenc}
\usepackage[russian,english,british]{babel}

\begin{document}

\title{Notation for Iteration of Functions, Iteral}
\author{Valerii Salov}
\date{}
\maketitle

\begin{abstract}
A new mathematical notation is proposed for the iteration of functions. It facilitates the application of the iteration of functions in mathematical and logical expressions, definitions of sets, and formulations of algorithms. Illustrations of the notation include definitions of constant points, periodic points, a filled-in Julia set, the Mandelbrot set, iterations of a logistic map, the double-approximating procedure for solving the Lorenz equations, a description of a financial time series, and reordering nonnegative integers useful for the investigation of the Collatz's (3x+1)/2 convergence problem. The terms iteral and iteral of function are suggested to name the new denomination.
\end{abstract}

\section{Introduction}

A mathematician or a natural science researcher routinely writing in equations $f(x), n!, \sum, \prod,$ and $\int$ relies on clarity, expressiveness, and convenience of mathematical notation filtered out by centuries of work, recognition, and acceptance. Creating a theory in many respects is developing a language eventually accompanied by a system of signs and symbols. A historical excursus \cite[pp. 441 - 443]{boyer} and \cite{cajori} shows that the road of many habitual mathematical symbols was bumpy with periods of oblivion, sometimes rediscovery, and final failure or success. A lack of a well understood notation creates a barrier for a wider application of a subject. In attempt to reuse the iteration of functions for the analysis of the financial high frequency time series the author did not find a notation suitably expressing this concept for writing formulas and equations.

Latin \textit{iterum} means \textit{again}. It is the root for the word \textit{iteration}. Iteration is a fundamental notion in computer science, where it denotes \textit{repetition of steps}, the repetition of a sequence of instructions in a computer program or algorithm. Accordingly, all computer languages have the means to express iterations. For instance, in C++ \cite[pp. 136 - 137]{stroustrup} these are \texttt{for}, \texttt{while}, and \texttt{do-while} iteration statements. The Unified Modeling Language abstracts the \textit{iterative expression} and places the iteration above concrete programming languages \cite[pp. 318 - 319]{rumbaugh}.

Dynamical systems, theory of chaos, and fractals depend on iterations. Benoit Mandelbrot comments \cite[p. 260, p. 268]{mandelbrot}:

\textit{The motivation for iteration and many of its important properties come straight from physics. In iteration, physicists see a simplified view of the state of a dynamical system and its evolution in time. \ldots The theory of iteration of rational functions goes back to the mid-nineteenth century and perhaps even to Abel.}

Marek Kuczma et al. \cite[p. 13]{kuczma} emphasize the relationship between iterations and functional equations:

\textit{We feel that any attempt to divorce iteration from functional equation investigations would be an extremely, indeed totally, fruitless task.}

This article proposes a new notation and term for iteration of functions and illustrates their usefulness.

\section{Existing Notations Implicitly Using Iterations}

The sign $\sum$ for summation was introduced by Leonard Euler in \textit{Institutiones calculi differentialis}, St. Petersburg, 1755, Cap. I, \$ 26, p. 27, who was saying, \textit {summam indicabimus signo} $\sum$ \cite[item 438, p. 61]{cajori}. In its variations the notation assigns initial values to the involved quantities, defines the number of iterations, binds each iteration with an arithmetic addition of the current term to a previous accumulated result, and assigns a value to the operation. It is unambiguous. What is going on in the following left-hand sides is intuitively well understood and explained by the middle expressions
\begin{displaymath}
\sum_{i=1}^{n} i = 1 + 2 + \cdots + n = \frac{n(n + 1)}{2},
\end{displaymath}
\begin{displaymath}
\sum_{i=1}^{\infty} \frac{1}{i^2} = \frac{1}{1} + \frac{1}{2^2} + \frac{1}{3^2} + \cdots = \frac{\pi^2}{6}.
\end{displaymath}
The expressions assume values, right-hand sides, and can become parts of other valid mathematical expressions. Presence of the counting index in the expression under the sum is not essential. Occasionally, the number of iterations is equal to zero. Then, a default, often zero, value is set for the sum.

Another successful sign for \textit{n}-factorial of natural numbers denoting "$n(n - 1)(n - 2) .... 3.2.1$" was introduced by Christian Kramp of Strasburg in his \textit{\'{E}l\'{e}mens d'arithm\'{e}tique universelle}, Colonge, 1808 \cite[item 448, p. 72]{cajori}. By convention $0! = 1$. The finite number of iterations is bound to the multiplication of factors incrementing by one and starting from one. The left-hand side below is more convenient than the middle or right-hand side
\begin{displaymath}
3! = 1 \times 2 \times 3 = \sum_{i=1}^{1}1 \times \sum_{i=1}^{2}1 \times \sum_{i=1}^{3}1.
\end{displaymath}
It can shorten descriptions: $29!$ vs. 8841761993739701954543616000000.

The factors in $n!$ are fixed for natural $n$. The signs for the product of $n$ arbitrary factors were introduced by Thomas Jarrett, who extensively studied algebraic notations and published results in 1831 \cite[item 448, pp. 69 - 71]{cajori}. His simplest proposal commented as "the product of $n$ factors, of which the $m$th is $a_m$" was $P_{m}^{n}a_{m}$. The $n$ was exactly above $P$. Today we write
\begin{displaymath}
\prod_{m=1}^{n} a_m.
\end{displaymath}
In 1812 Carl Friedrich Gauss introduced $\prod$ in the designation of continued infinite products \cite[item 451, p. 78]{cajori}. The symbols $\Pi{z}$ were used by Gauss in \textit{Werke}, Vol. III, G\"{o}ttingen, 1866, p. 146 for expressing the gamma function \cite[item 650, p. 272]{cajori}. Camill Jordan applies the notation $\prod_r(c^{\pi{p^r}}-K^\pi)^e$ for the multiplication of factors resulting from different values of $r$ published in 1870. George Chrystal represents the product $(1 + u_1)(1 + u_2)....(1 + u_n)$ by $\prod^n{1 + u_n}$ printed in 1889 \cite[item 451, p. 79]{cajori}. The home page of the Clay Mathematics Institute displays a photocopy of the first page of the hand written manuscript of Bernhard Riemann's paper \textit{Ueber die Anzahl der Primzahlen unter einer gegebenen Gr\"{o}sse} dated by 1859, November and first time published in \textit{Monatsberichte der Berliner Akademie} in the same month \cite{clay}. The uppercase cursive Greek letters $\Pi$ and $\Sigma$ used without indexes build the starting Riemann's point, the identity proved by Euler and known as the Euler product formula
\begin{displaymath}
\prod{\frac{1}{1 - \frac{1}{p^s}}} = \sum{\frac{1}{n^s}}
\end{displaymath}
Here $p$ ranges over all primes 2, 3, 5, 7, ... and $n$ ranges over all positive integers 1, 2, 3, 4, ..., what is frequently symbolized by placing $p$ and $n$ under the product and sum signs. The right-hand side contains infinite number of terms and diverges for $s = 1$. Since the left-hand side involves unique primes without repetitions it would converge for any finite number of factors. This implies infinite number of primes. Jay Goldman comments \cite[p. 36]{goldman}:

\textit{This may seem a rather difficult way of proving a theorem that Euclid had already proved by very elementary reasoning, but in the 19th century Dirichlet vastly extended the ideas used here and they form the basis for the most important applications of analysis in modern number theory.}

For us it means that Riemann similarly to Gauss applies the product sign for the infinite number of factors. In \cite{boyer}, \cite{cajori}, and \cite{yushkevich} the author did not find, when and by whom the symbol $\prod$ was used first time for the product of $n$ factors. However, the popular mathematical handbook of the second half of the 20th century lists it among the first mathematical denominations in introductory sections \cite{korn}. The product with zero iterations is assigned the default value one. 

A well known notation combining the limit with the sum of terms is $\int$. Most frequently it implies infinite number of iterations. Jacob Bernoulli was the first, who used the word \textit{integral} in print. Johann Bernoulli claimed for himself the introduction of the word. The principal symbol of the integral calculus has been suggested by Gottfried Leibniz as the long letter $\int$ reminding about \textit{calculus summatorius}, the name favored by Leibniz over the name \textit{calculus integralis} and the capital letter \textit{I} favored by Johann \cite[item 539, pp. 181 - 182]{cajori}. The sign $\int$ is successfully applied to denote deterministic Riemann, Stieltjes, Lebesgue, Riemann-Stieltjes, Riemann-Lebesgue \cite{korn}, and stochastic It\^{o}, and Stratonovich integrals \cite[p. 48, p. 106]{rogers}. Frequently, in abstract formulations under the integral sign we find useful for us combination of characters $f(x)$.

The first appearance of $f$ for "function" together with parentheses is dated by 1734 and reported for Euler's \textit{Comment. Petropol. ad annos 1734-1735}, Vol. VII, 1840, p. 186. This is cited by a secondary reference of Johannes Tropfke in \cite[item 643, p. 268]{cajori} as the Euler's words: "Si $f(\frac{x}{a}+c)$ denotet functionem quamcunque ipsius $\frac{x}{a}+c$". Let us list the following common and useful properties of the discussed notations

\begin{itemize}
\item The symbols $\sum, n!, \prod, \int$ imply finite or infinite iterations.
\item The number of iterations and initial quantities are defined.
\item The $f(x)$ under $\int$ is visible and can be substituted by a concrete function.
\item The denominations assume values and can be parts of other expressions.
\item The notations are unambiguous and can be mixed.
\end{itemize}

What prevents applying the symbols based on iterations for defining the iteration of functions is binding $\sum$ and $\int$ to addition and $\prod$ and $n!$ to multiplication. Let us review existing notations for the iteration of functions.

\section{Existing Notations for Iterations of Functions}

Articles and monographs are written on functional equations, dynamical systems, chaos, fractals, and programming languages relying on iterations. The operation of composition is the basis of the key definitions. According to \cite{kuczma}, this is the only inner operation which can be defined in the family of all self-mappings $\mathcal{F}(X)$ of a given set $X$. \textit{All} also includes such an element $f \in \mathcal{F}(X)$, denoted $\mathrm{id}_X$, that $f(x) = x$ $\forall x \in X$. Since the composition $\circ$ is associative and $f \circ \mathrm{id}_X = \mathrm{id}_X \circ f$ $\forall f \in \mathcal{F}(X)$, the system $(\mathcal{F}(X), \circ)$ forms a semigroup with identity $\mathrm{id}_X$, monoid. In general, it is noncommutative. With these properties the \textit{iterates} are defined as the powers $f^n$ of an element $f \in \mathcal{F}(X)$, where $n$ is a positive or zero integer
\begin{displaymath}
n \in \mathrm{N}_0, f^0 = \mathrm{id}_X, f^{n + 1} = f \circ f^n.
\end{displaymath}
While two arbitrary elements of the family do not have to commute with respect to the composition, the iterates of a function $f \in \mathcal{F}(X)$ do
\begin{displaymath}
n, m \in \mathrm{N}_0, f^n \circ f^m = f^m \circ f^n = f^{n + m}.
\end{displaymath}
Here we are switching from the terms \textit{element} and \textit{self-mapping} to \textit{function} using Euler's notation $f(x)$. The composition of two functions $f(x)$ and $g(x)$ is expressed as $(f \circ g) = f(g(x))$.

The $x$ is a \textit{predecessor} of $y$ under $f$ and $y$ is a \textit{successor} of $x$ under $f$, if there exists an $n \in \mathrm{N}_0$ such that $y = f^n(x)$. The ordering between predecessor and successor is partial. It is reflexive, and transitive but not antisymmetric. For example, if $x \neq 1$ and $f(x) = \frac{1}{x}$, then $b = f(a), a = f(b),$ and $a \neq b$. The $x$ and $y$ are \textit{equivalent}, if they have a common successor under $f$: $f^n(x) = f^m(y)$. This is an \textit{equivalence relation}, where reflexivity and symmetry are obvious and transitivity follows from $f^{m + k}(x) = f^{m + l}(y) = f^{l + m}(y) = f^{l + n}(z)$. The equivalence classes under $f$ are called \textit{orbits} \cite[pp. 13 - 18]{targonski}. This relation is fundamental for the theory of iterations. Gy\"{o}rgy Targonski \cite[p. 14, p. 283]{targonski} describes that it was introduced in a short remark made by Kazimierz Kuratowski in the text of a paper written in 1924 by Ralph Tambe-Lyche. It was independently rediscovered 17 years later and printed as a 14 lines abstract in 1941 by Gordon Thomas Whyburn, who coined the terms \textit{orbit} and \textit{orbit decompositions} and began systematic study of the new equivalence classes \cite[p. 14]{kuczma}. Possibility to select $x$, keep it intact, and change the power $n \in \mathrm{N}_0$ sequentially creates the set $\{f^n(x)\}$ called the \textit{iteration sequence}. Targonski \cite{targonski} suggests to name it \textit{splinter} of $x$ borrowing the term from the theory of recursive functions, where it was introduced by Joseph Ullian \cite{ullian}.

The book \cite{targonski} is considered in \cite[p. 14]{kuczma} as the first monograph completely devoted to the iteration theory. Systematically applying $f^n(x)$ and $\circ$ it successfully develops the theory proving notation's sufficiency. In practical applications the set $X$ usually has an algebraic structure, where frequently it is needed to distinguish the iterates and powers with respect to multiplication. Anticipating such a complication the authors of \cite{kuczma} propose $f^{2}(x)$ for the composition $f[f(x)]$ and $[f(x)]^2$ for the power $f(x)f(x)$. Something is needed. Otherwise, what is the value of $\cos^2{x} + \sin^2{x}$ in a mixed environment? It is unity for powers and something else, if at least one term assumes composition. Following the proposal the mixed expression $\cos^2{x} + [\sin{x}]^2$ would be confusing because the parentheses are general grouping symbols. What have others applied in the historical order of publications?

Olexandr (Alexander) Sharkovsky deriving properties of a continuous mapping $T$ of the line into itself denotes iterates $T^k{\alpha}$ \cite{sharkovsky1}. He switches to lower and upper indexes for an infinite sequence of iterated values of $x$, $\{T^j{x}\}_{j = 0}^\infty$ \cite{sharkovsky2}.

Tien-Yien Li and James Yorke \cite{li} write $x_{n + 1} = F(x_n)$ and expand it for the logistic equation $x_{n + 1} = rx_n[1 - x_n/K]$. One recursive step is represented well. However, to use it in expressions, where initial values are assumed and the number of iterations is arbitrary is problematic.

Mitchell Feigenbaum contributing to nonlinear dynamics applies the same $x_{n+1} = f(x_n)$ for recursion equations \cite{feigenbaum}. The function $f(x)$ can be easily replaced by a concrete one. However, substituting either the left-hand or right-hand side into other expressions either can be interpreted as vector coordinates or loses the number of iterations, when it is greater than one.

Robert Devaney presents several programs for visualization of notions of chaos, fractals, and dynamical systems. Computer science demands iterations. Devaney exploits two step definitions like: for $x > 1$ and $T(x) = x^2, T^n(x) \to \infty$ as $n \to \infty$ \cite[p. 12]{devaney}. For instance, $T^{3}(1.2)$ gets the value 4.29981696. It is not clear how to replace $T$ by $x^2$. A reader can  either be puzzled by $(x^2)^3(1.2)$ or give the wrong answer 2.985984 for $(1.2^2)^3$. Doing the composition in two steps is similar to writing $P(x) = a + bx + cx^2$ and $\int{P(x)dx}$ instead of $\int{(a + bx + cx^2)dx}$.

Donald Knuth working with iteration of series favors $x_n \gets f(x_{n - 1})$ for an iterative process and $x_n = f(f(...f(x_0)...))$ for the $n$-fold composition of a given function $f$ \cite[p. 530]{knuth}. He chooses square parentheses $f^{[n]}(x) = f(f^{[n - 1]}(x))$, reducing ambiguity with powers of multiplication. Sequential values $x_n$ can be misinterpreted as vector coordinates. Replacing $f^{[x]}(x)$ by an actual function can be misused. For the polynomial from the previous paragraph we get $P^{[n]}(x) = (a + bx + cx^2)^{[n]}$ what is not intuitive.

Mandelbrot considers the map $z \to \lambda z(1 - z)$ and iterated maps $z_n = f(f(....f(z_0)))$ \cite[p. 38]{mandelbrot}. The last right-hand side is good for illustrations of iterates and not suitable for expressions. It is not concrete for referencing the number of iterations. Using the concrete left-hand side $z_n$ would be ambiguous in cases involving vector coordinates.

We see that there is no common notation for the iteration of functions. Existing ones address theoretical needs and become ambiguous or less intuitive in mixed application environments.

\section{Proposal of Iteral}

When a mathematical notion is getting fundamental and other disciplines are demanding it, it may be a time for a common notation. The iteration of functions deserves it. The following are the simplest proposals
\begin{displaymath}
\mathrm{SYMBOL}_\mathrm{INITIAL \; VALUE}^\mathrm{ITERATIONS \; PERFORMED}f(x) \; \textrm{or}
\end{displaymath}
\begin{displaymath}
\mathrm{SYMBOL}_{x=\mathrm{INITIAL \; VALUE}}^\mathrm{ITERATIONS \; PERFORMED}f(x)
\end{displaymath}

Here, $\mathrm{SYMBOL}$ should correspond to the meaning of iteration, $\mathrm{INITIAL}$ $\mathrm{VALUE}$ is an element of the domain set of the function $f(x)$, $\mathrm{ITERATIONS}$ $\mathrm{PERFORMED}$ is the number of iterations completed so far and started from 0.

If $\mathrm{ITERATIONS \; PERFORMED} = 0$, then the value of the expression $= \mathrm{INITIAL \; VALUE}$. A value of the function $f(x)$ after a previous iteration, or $ \mathrm{INITIAL \; VALUE}$ after the beginning, is a new argument of the function for the next iteration. Iterations may continue up to $\mathrm{ITERATIONS \; PERFORMED}$ only if intermediate values are members of the domain set of the function $f(x)$. The second version explicitly indicates the variable $x$, where and expression for $f(x)$ makes it unclear.  This notation combines the following useful properties:

\begin{itemize}
\item The function is visible and can be replaced with a concrete one.
\item The initial value is visible and can be replaced with a concrete expression.
\item The number of iterations is unambiguously specified.
\item The expression has a value and can be a part of other expressions.
\item Combination of structural parts of $\sum$ and $\int f(x)dx$ is intuitive.
\end{itemize}

Greek letters are frequent candidates for symbols. To be unambiguous, they should differ from English letters. By this reason, Greek upper case iota is not suitable. Invention of a new symbol like $\int$ is challenging. The last reminds about \textit{summatorius} and our should remind about \textit{iterum} and \textit{iteration}. Latin \textrm{\"I} has distinguishing dots but it is not as exact as \textrm{I}. Let us create a mirror image of \textrm{N} $\to$ \CYRI. This makes it the first Cyrillic letter in Russian \textit{\CYRI\cyrt\cyre\cyrr\cyra\cyrc\cyri\cyrya}, translation of English \textit{Iteration}, and a reminder. Thus,
\begin{displaymath}
\textrm{\LARGE{\CYRI}}_{v}^{n}f(x) \; \textrm{or} \; \textrm{\LARGE{\CYRI}}_{x=v}^{n}f(x)
\end{displaymath}
are the new notations. The author suggests to name this operation the \textit{iteral of function} or simply \textit{iteral}. This correlates with the formation rules of the words \textit{integral, differential, factorial,} and \textit{fractal}. Due to the Latin root, \textit{iteral} can be naturally consumed by other languages.

\section{Applications of Iterals}

Below $(\ldots)$ and $\{\ldots\}$ denote a sequence and a set. $\mathrm{Z, N, R, C}$ are sets of integers, positive integers, real, and complex numbers. $\mathrm{N_0} = \{0, \mathrm{N}\}$. $(\mathrm{N_0}) = (0, 1, 2, ... )$ is the infinite sequence. $(n \in \mathrm{N_0}) = (0, 1, 2, ..., n)$ is a finite sequence, \textit{chain}, of nonnegative integers. The \textit{set builder notation} $\{x : P(x)\}$ means the set of all $x$ for which $P(x)$ is true. The \textit{universal quantification of proposition} $\forall$ "for all", \textit{the existential quantification of proposition} $\exists$ "there exists", $\in$ "member of", $\wedge$ "and", $\vee$ "or" are standard. 

\textbf{Iterals have values.} To get flexibility: \CYRI$_2^0{x^2} = 2$, \CYRI$_2^1{x^2} = 4$, \CYRI$_2^2{x^2} = 16$, \CYRI$_2^\infty{x^2} \to \infty$, \CYRI$_1^\infty{x^2} = 1$, \CYRI$_{0.5}^\infty{x^2} \to 0$. For $x = \frac{\pi}{4}, \cos^2{x} + $ \CYRI$_{x}^2\sin x = \frac{1}{2} + \sin \frac{\sqrt{2}}{2}$. Here are more iterals
\begin{displaymath}
\textrm{\LARGE{\CYRI}}_{0}^{n}(x + 1) = n, \textrm{\LARGE{\CYRI}}_{x=1}^{n}ax = a^n, \textrm{\LARGE{\CYRI}}_{1}^{\infty}\frac{1}{x + 1} = \lim_{n \to \infty} \textrm{\LARGE{\CYRI}}_{1}^{n}\frac{1}{x + 1} = \frac{2}{1 + \sqrt{5}}
\end{displaymath}
The last is reciprocal of the \textit{golden ratio}. Indeed, $\forall (n \in \mathrm{N_0})$ the iterals form the chain $(\frac{1}{1}, \frac{1}{2}, \frac{2}{3}, \frac{3}{5}, \frac{5}{8}, \frac{8}{13}, \; \ldots)$, where each member is the ratio of a previous Fibonacci's number to the current one. Using Targonski's and Ullian's term proposal, the splinter of $1$ for the function $\frac{1}{x + 1}$ is the set $\{$\CYRI$_{1}^{n}\frac{1}{x + 1}\} \; (n \in \mathrm{N_0})$. The splinter is an infinite sequence, \textit{acyclic splinter}, if $\forall m,n \in \mathrm{N} \; m \neq n$ \CYRI$_{x}^m f(x) \neq$ \CYRI$_{x}^n f(x)$. The splinter is cyclic, if $\exists m,n \in \mathrm{N} \; m < n$ \CYRI$_{x}^m f(x) =$ \CYRI$_{x}^n f(x)$. 

\textbf{Iterals form algebraic expressions and functions.} Since the initial value of the iteral substitutes the argument of the function under the iteral sign the notation creates algebraic expressions or functions, if the initial value contains parameters and/or the function argument. The composition of functions must be valid. Iterals can nest. Here are a few examples used later for reordering $\mathrm{N_0}$.
\begin{displaymath}
\forall n \in \mathrm{N_0} \; \textrm{\LARGE{\CYRI}}_{v}^{n}p = v = \mathrm{id}_v, \;  \textrm{\LARGE{\CYRI}}_{2p}^{n}p = 2p, \; \textrm{\LARGE{\CYRI}}_{2p+1}^{n}p = 2p+1,
\end{displaymath}
\begin{displaymath}
\forall k \in \mathrm{N_0} \; \textrm{\LARGE{\CYRI}}_{2p+1}^{k}(2p+1) = 2^{k+1}p+2^{k+1}-1,
\end{displaymath}
The last is proved by induction. For $k = 0$ it holds: \CYRI$_{2p+1}^{0}(2p+1) = 2p + 1 = 2^{0+1}p+2^{0+1}-1$. If it holds for $k$, then for $k + 1$ the left-hand side is equal to
\begin{displaymath}
\textrm{\LARGE{\CYRI}}_{2p+1}^{k+1}(2p+1) = \textrm{\LARGE{\CYRI}}_{2^{k+1}p+2^{k+1}-1}^{1}(2p+1) = 2^{k+2}p+2^{k+2}-1,
\end{displaymath}
and the right hand side is also equal to $2^{(k + 1)+1}p+2^{(k + 1)+1}-1 = 2^{k + 2}p+2^{k + 2}-1$.

The proved iterals describe \textit{unfixed subsequences} with the names $O_k{O}$ defined in the section on subdivision of $\mathrm{N_0}$. The index $k$ is the number of prefixing extra characters $O$, meaning \textit{odd}, in the name of odd numbers $O$. For $k = 0$ the name is just $O$, and the subsequence, given by the formula $2p+1$, represents all positive odd numbers with positions $p = 0, 1, \; \ldots$ within this subsequence. The following will describe the \textit{fixed subsequences} $EO_k{O}$ of the odd numbers
\begin{displaymath}
\textrm{\LARGE{\CYRI}}_{2p}^{1}
\textrm{\LARGE{\CYRI}}_{2p+1}^{k}(2p+1) = 2^{k+2}p+2^{k+1}-1, \; k = 0, 1, \; \ldots \; .
\end{displaymath}
The $E$ associates with \textit{even}. The fixed subsequence of odd numbers with the shortest name $EO$ for $k = 0$ is defined by $4p+1$. Geometrically, the lines $2^{k+2}p+2^{k+1}-1$ of different orders $k$ do not intersect for $p \geq 0$. 

In general, nesting or composing iterals for the same function with common initial values is not the same as increasing the number of iterations
\begin{displaymath}
\textrm{\LARGE{\CYRI}}_{2p+1}^{1}\textrm{\LARGE{\CYRI}}_{2p+1}^{1}2p =
\textrm{\LARGE{\CYRI}}_{2p+1}^{1}(4p+2) = 8p + 6,
\end{displaymath}
\begin{displaymath}
\textrm{\LARGE{\CYRI}}_{2p+1}^{2}2p = \textrm{\LARGE{\CYRI}}_{4p+2}^{1}2p = 8p + 4.
\end{displaymath}
This is because going from iteration to iteration the iteral's function remains intact and the initial value changes
\begin{displaymath}
\textrm{\LARGE{\CYRI}}_{v}^{k+1}f(x) = \textrm{\LARGE{\CYRI}}_{\textrm{\LARGE{\CYRI}}_{v}^{k}f(x)}^{1}f(x).
\end{displaymath}
While nested iterals keep the common initial value intact and likely change the function on each composition step.

We need more examples for the section on subdivision of $\mathrm{N_0}$
\begin{displaymath}
\forall m \in \mathrm{N_0} \; \textrm{\LARGE{\CYRI}}_{2p}^{m}2p = 2^{m+1}p.
\end{displaymath}
For $m = 0$ it holds by the definition of iteral \CYRI$_{2p}^{0}2p = 2p$ and $2^{0+1}p = 2p$. The induction step from $m \to m + 1$ gives the same result for the left- and right-hand sides: \CYRI$_{2p}^{m+1}2p = $ \CYRI$_{2^{m + 1}p}^{1}2p =2^{m+2}p$ and $2^{(m+1)+1}p = 2^{m+2}p$.
\begin{displaymath}
\forall m, l \in \mathrm{N_0} \; \left(\textrm{\LARGE{\CYRI}}_{2p+1}^{1}\cdots
\textrm{\LARGE{\CYRI}}_{2p+1}^{1}\right)_{l}
\textrm{\LARGE{\CYRI}}_{2p}^{m}2p = 2^{m + l + 1}p + 2^{m+1}(2^l -1),
\end{displaymath}
The index $l$ is the number of composing iterals with the common initial value $2p + 1$. Due to the previous proof the left-hand side is equal to
\begin{displaymath}
\left(\textrm{\LARGE{\CYRI}}_{2p+1}^{1}\cdots
\textrm{\LARGE{\CYRI}}_{2p+1}^{1}\right)_{l}
2^{m + 1}p.
\end{displaymath}
For $l = 0$ the equality holds. In this case the iterals in the left-hand side vanish leaving $2^{m + 1}p = 2^{m + 0 + 1}p + 2^{m+1}(2^0 -1) = 2^{m + 1}p + 0 = 2^{m + 1}p$. The equality holds for $l = 1$ 
\begin{displaymath}
\textrm{\LARGE{\CYRI}}_{2p+1}^{1}2^{m + 1}p = 2^{m + 2}p + 2^{m + 1} \; \textrm{and} \;
2^{m + 1 + 1}p + 2^{m+1}(2^1 -1) = 2^{m + 2}p + 2^{m + 1}
\end{displaymath}
On the induction step $l \to l + 1$ the left-hand side gives
\begin{displaymath}
\left(\textrm{\LARGE{\CYRI}}_{2p+1}^{1}\cdots
\textrm{\LARGE{\CYRI}}_{2p+1}^{1}\right)_{l+1}
2^{m + 1}p = \textrm{\LARGE{\CYRI}}_{2p+1}^{1}(2^{m + l + 1}p + 2^{m+1}(2^l -1)) =
\end{displaymath}
\begin{displaymath}
= 2^{m + l + 2}p + 2^{m + l + 1} + 2^{m+1}(2^l -1) = 2^{m + l + 2}p + 2^{m+1}(2^l + 2^l -1) =
\end{displaymath}
\begin{displaymath}
= 2^{m + l + 2}p + 2^{m+1}(2^{l+1} -1).
\end{displaymath}
The same result is for the right-hand side $2^{m + (l + 1) + 1}p + 2^{m+1}(2^{(l+1)} -1) = $
$2^{m + l + 2}p + 2^{m+1}(2^{l+1} -1)$. Then, the following result is simple
\begin{displaymath}
\forall m, l \in \mathrm{N_0} \; \textrm{\LARGE{\CYRI}}_{2p}^{1}\left(\textrm{\LARGE{\CYRI}}_{2p+1}^{1}\cdots
\textrm{\LARGE{\CYRI}}_{2p+1}^{1}\right)_{l}
\textrm{\LARGE{\CYRI}}_{2p}^{m}2p = 2^{m + l + 2}p + 2^{m+1}(2^l -1),
\end{displaymath}
For $l \geq 1$ the last iterals correspond to fixed subsequences of even numbers defined later and named $EO_l{E}_m{E}$. The shortest name $EOE$ for $l = 1$ and $m = 0$ corresponds to the sequence $8p + 2$.

Similar to tables of sums, products, and integrals it would be useful to have tables of iterals and add iterals into symbolic processing computer systems.

\textbf{Iterals work in definitions and statements}. \textit{Fixed point $v$}: $f(v) = v$. Using iterals: \CYRI$_{v}^{1}f(x) = v$. The original is shorter. However, the following property of a fixed point is expressed unambiguously with iterals: $\forall n \in \mathrm{N}$ \CYRI$_{v}^n f(x) = v$. Examples: $\forall n \in \mathrm{N}$ \CYRI$_{1}^n x^2 = 1$ and \CYRI$_{0}^n x^2 = 0$.

\textit{Periodic point $v$ of the order $m$}: $\exists k,m  \in \mathrm{N} \; \forall k<m$ \CYRI$_{v}^m f(x) = v \; \wedge$ \CYRI$_{v}^k f(x) \neq v$. Example: $\forall v \in R \wedge v \neq 0$ \CYRI$_{v}^2\frac{1}{x}  = v \; \wedge$ \CYRI$_{v}^1\frac{1}{x}  \neq v$.

\textit{A filled in Julia set}. Biographical facts of Pierre Fatou and Gaston Julia and historical comments on iterations of rational functions, particularly, the complex function $f(z) = z^2 + c$, can be found in \cite[pp. 268 - 275]{mandelbrot}. The iteral definition of the set is: $z \in \mathrm{C}, \; J_c = \{v \in \mathrm{C} : \; $|\CYRI$_{v}^{\infty}(z^2 + c)| < \infty\}$.

\textit{The Mandelbrot set}. $z \in \mathrm{C}, \; M = \{c \in \mathrm{C} : \; $|\CYRI$_{0}^{\infty}(z^2 + c)| < \infty\}$. The iteral definitions of the last two sets give straight forward directions for writing computer programs depicting the points of the complex plane. C++ provides for the class \texttt{complex} making arithmetic and functional programming involving the complex numbers natural \cite[pp. 267 - 274]{stroustrup}.

\textit{Iterations of logistic map}. Two Belgian mathematicians Lamberte Adoplhe Jacques Quetelet and Pierre Fran\c{c}ois Verhulst contributed to modeling of the growth of a biological population using an ordinary differential equation with separable variables of the population size and time. A solution of it is an S-shaped function. Ka-Kit Tung \cite[p. 100]{tung} writes:

\textit{Verhulst called this solution the "logistic curve". The French term "logistique" was used to signify the art of calculation.}

In a form of recursive equation it is named the \textit{logistic map} \cite{feigenbaum} and was studied in \cite{li}, where Yorke has coined the term $chaos$. According to the authors, chaos is a result of presence of periodic orbits with all periods. Getting spectacular results the authors did not know about Sharkovsky's theorem \cite{misiurewicz}. Yakov Sinai \cite{sinai} marks historical milestones of the chaos theory and anticipates its future:

\textit{I believe that the future of the chaos theory will be connected with new phenomena in non-linear PDE and other infinite-dimensional dynamical systems, where we can encounter absolutely unexpected phenomena.}

This development existed long before the introduction of the term chaos. Vladimir Arnold \cite{arnold} emphasizes the origins in the work of Henry Poincar\'{e}. If we narrow the topic, then the orbits of a logistic map are given by \CYRI$_{v}^{n}bx(1-x)$. A chaotic time series is extremely sensitive to the initial conditions. Rounding off errors of computer operations evaluating initial conditions before a next iteration can become a natural source of instability in modeling such systems.

\textit{Lorenz equations.} Edward Lorenz \cite{lorenz}, using a modified model of Barry Saltzman \cite{saltzman} for simulation of atmospheric phenomenon, have discovered unexpected behavior of a solution of the system of three ordinary differential equations referred today as the Lorenz equations
\begin{displaymath}
\frac{dX}{d\tau} = -\sigma X + \sigma Y, \; \frac{dY}{d\tau} = -XZ + rX - Y, \; \frac{dZ}{d\tau} = XY - bZ,
\end{displaymath}
where $\sigma, r,$ and $b$ are positive constants, and $\tau$ is a dimensionless time. The solution is a phase trajectory $(X(\tau), Y(\tau), Z(\tau))$, where $X, Y,$ and $Z$ define a phase space of a layer of fluid of uniform depth between two surfaces maintained at two different temperatures. Depending on the conditions and \textit{Rayleigh number} the liquid remains steady or gets in motion, \textit{convection}. $X$ is proportional to the intensity of the convection, $Y$ is proportional to the temperature difference between the ascending and descending currents, and $Z$ is proportional to distortion of the vertical temperature profile from linearity. While the equations can be solved using different numerical methods, the author reviewed the \textit{double-approximation procedure} described by Lorenz and found its iteral expression. Consider the initial value $\boldsymbol{P_0}$ of the vector variable $\boldsymbol{P}$, and another variable $\boldsymbol{p}$
\begin{displaymath}
\boldsymbol{P_0} =
\left(
\begin{array}{r}
X_0\\
Y_0\\
Z_0
\end{array}
\right); \;
\boldsymbol{P} =
\left(
\begin{array}{r}
X\\
Y\\
Z
\end{array}
\right); \;
\boldsymbol{p} =
\left(
\begin{array}{r}
x\\
y\\
z
\end{array}
\right).
\end{displaymath}
Then, a point of the \textit{approximated} trajectory after the $n$ iterations is equal to
\begin{displaymath}
\boldsymbol{P_n} = \textrm{\LARGE{\CYRI}}_{\boldsymbol{P_0}}^{n}0.5
(
\left(
\begin{array}{r}
X\\
Y\\
Z
\end{array}
\right) + \textrm{\LARGE{\CYRI}}_{\boldsymbol{P}}^{2}(
\left(
\begin{array}{rrr}
-\sigma&\sigma&0\\
r&-1&-x\\
0&x&-b
\end{array}
\right)
\left(
\begin{array}{r}
x\\
y\\
z
\end{array}
\right)\Delta\tau +
\left(
\begin{array}{r}
x\\
y\\
z
\end{array}
\right)
)).
\end{displaymath}
This is a direct way to a computer program.

\textit{A financial time series.} The most complete information about transactions comes from the markets as ticks, triplets of time, price, and number of contracts or shares named volume $(t_i, P_i, V_i)$. They are ordered by time. An ordinary daily session of a liquid nearby S\&P500 E-mini futures contract traded electronically on GLOBEX platform on the Chicago Mercantile Exchange brings several hundred thousands ticks. This associates with the name \textit{high-frequency trading} \cite{engle}. These data is a result of cooperation of modern technology and human consciousness governed by partly unknown laws of nature. Substantial and frequent potential profits attract to trading. Real losses sober. Inability to avoid trading risk and the behavior of the market "attempting" to fool the majority of its participants creates a big assortment of approaches for "breaking the market code" spanning from science to astrology. The markets appear random. A hope of traders is to find something more than a pure luck but a dependence between past and current behavior. Since all cannot win together, because the market does not create treasure but redistributes it giving away a portion to the industry supporting trading, \textit{non-zero sum game}, the task of systematic winning by minor percent of participants becomes extremely tough.

Probability theory, statistics, and theory of stochastic processes are traditional ways to study the markets. Behavioral finance attempts to uncover psychological contribution into the process and find suitable quantitative measures of it \cite{tversky}. Forgetting about human being component of the activity can be a costly trading experiment. Albert Shiryaev calls the \textit{first order task} studying statistics of waiting times, $\Delta t_i = t_{i} - t_{i-1}$ together with price increments, $\Delta P_i = P_{i} - P_{i-1}$ \cite[p. 379 of Russian edition]{shiryaev}. More and more empirical observations confirm that Bachelier's or Samuelson's representations of price increments or logarithms of price ratios, asset returns, by random Gaussian variables are only simplifications \cite{salov}. Robert Engle describes \cite{engle} that the Weibull \cite{weibull} distribution well approximates empirical waiting times. This as well as the generalized Kumaraswamy distribution \cite{kumaraswamy} are suitable for the waiting times coming from the Chicago Board of Trade \cite{salov}. Discrete probability distributions and laws such as multinomial, Hurwitz Zeta, and Zipf-Mandelbrot can approximate discrete price increments and their extreme values \cite{salov}. Prices and their increments are essentially discrete and can be expressed as whole numbers of ticks. The author names the indecomposable further waiting times and price changes coming from a single trading session as \textit{a-} and \textit{b-increments} and the changes between the first, in a next session, and the last, in a current session, prices \textit{c-increments}. These are quantitative measures of the a-, b-, and c-properties. Thus, the a-b-c-process evolves as following: a) the first property determines the moment of the next transaction in a trading session; b) the second property carries on the first one and determines the price fluctuation expressed by a whole number of ticks; and c) the third property is responsible for the discrete price change between the current last and next first prices. Selecting time $u_t$ and price $u_p$ units for depicting an a-b-c-process on a chart, each tick becomes a point $(t_i{u_t}, p_i{u_p})$. If this is viewed as a vector or a complex number with imaginary price component, then $\Delta z_i = ((t_i - t_{i - 1})u_t, (p_i - p_{i - 1})u_p) = (u_t \Delta t_i, u_p \Delta p_i)$. Let $Z_0$ be the last tick of a previous session, $\Delta Z(\Delta A, \Delta C)$ be the random c-increment between a next first and a last previous ticks, and the lower case $z$ be responsible for the intraday ticks. Then, using the iterals the a-b-c process is
\begin{displaymath}
\textrm{\LARGE{\CYRI}}_{Z_0}^{j} \textrm{\LARGE{\CYRI}}_{Z + \Delta Z}^{N_j - 1}(z + \Delta z),  
\end{displaymath}
where $j = 0, 1, 2, \ldots$ enumerates sessions, $N_j$ is a random number of intraday ticks, and $\Delta z$ is a random variable. To which degree everything is random, independent, or varying in time is the subject of numerous modern investigations.
 
\textbf{A division of $\mathrm{N}_0$ on subsequences. A sieve.} Below are the sequences of the nonnegative integers $n$, their positions $p$, and even $e$ or odd $o$ types $t$ 
\begin{displaymath}
\begin{array}{rrrrrrrrrrrrrrrrrrr}
n & 0 & 1 & 2 & 3 & 4 & 5 & 6 & 7 & 8 & 9 & 10 & 11 & 12 & 13 & 14 & \ldots & p & \ldots\\
p & 0 & 1 & 2 & 3 & 4 & 5 & 6 & 7 & 8 & 9 & 10 & 11 & 12 & 13 & 14 & \ldots & p & \ldots\\
t & e & o & e & o & e & o & e & o & e & o & e & o & e & o & e & \ldots & e \vee o & \ldots
\end{array}
\end{displaymath}
A sieve is applied to the sequence $n$. The sieve extracts numbers occupying even positions and puts them into the subsequence $E$, where their positions are renumbered in the natural order starting from zero. The extended types indicate \textit{even-even ee} and \textit{odd-even oe} even numbers located on even and odd positions within $E$. The $E$ contains all and only even numbers including zero.
\begin{displaymath}
\begin{array}{rrrrrrrrrrrrr}
n & 0 & 2 & 4 & 6 & 8 & 10 & 12 & 14 & 16 & \ldots & 2p & \ldots \\
p & 0 & 1 & 2 & 3 & 4 &   5 &   6 &   7 &   8 & \ldots &   p & \ldots \\
t & ee & oe & ee & oe & ee & oe &  ee &   oe &  ee & \ldots &  ee \vee oe & \ldots
\end{array}
\end{displaymath}
The remaining numbers occupying odd positions in the original sequence $n$ are put into the subsequence $O$, where they get new positions counted also from zero. The $O$ contains all and only odd numbers. The odd numbers on even and odd positions within $O$ get \textit{even-odd eo} and \textit{odd-odd oo} extended types.
\begin{displaymath}
\begin{array}{rrrrrrrrrrrrr}
n & 1 & 3 & 5 & 7 & 9 & 11 & 13 & 15 & 17 & \ldots & 2p + 1 &\ldots \\
p & 0 & 1 & 2 & 3 & 4 &   5 &   6 &   7 &   8 & \ldots & p & \ldots \\
t &  eo & oo &  eo & oo & eo & oo & eo & oo & eo & \ldots & eo \vee oo & \ldots
\end{array}
\end{displaymath}
Being applied to $E$ the sieve splits it into $EE$ and $OE$. The new letter $E$ or $O$, depending on the evenness of the extracted positions, is concatenated to a previous name on the left side. Similarly, from $O$ the sieve produces $EO$ and $OO$ subsequences. The sieve is not applied to the subsequences getting the prefix $EO$ in their names. These subsequences are referred to as the \textit{fixed subsequences}. It is applied once to all other subsequences creating the following growing binary tree with stopping, final terminal, nodes corresponding to the fixed subsequences.
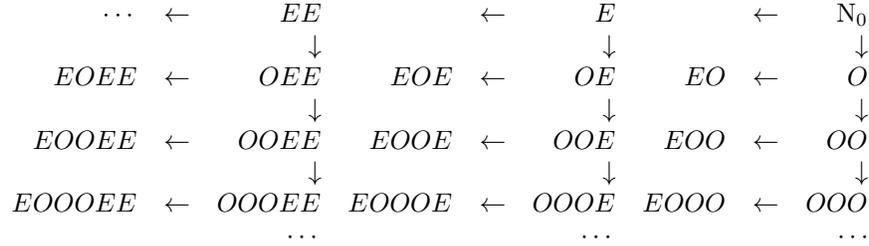
\begin{displaymath}
\begin{array}{rrrrrrrrr}
\cdots&\gets&EE&&\gets&E&&\gets&\mathrm{N_0}\\
&&\downarrow&&&\downarrow&&&\downarrow\\
EOEE&\gets&OEE&EOE&\gets&OE&EO&\gets&O\\
&&\downarrow&&&\downarrow&&&\downarrow\\
EOOEE&\gets&OOEE&EOOE&\gets&OOE&EOO&\gets&OO\\
&&\downarrow&&&\downarrow&&&\downarrow\\
EOOOEE&\gets&OOOEE&EOOOE&\gets&OOOE&EOOO&\gets&OOO\\
&&\cdots&&&\cdots&&&\cdots
\end{array}
\end{displaymath}
To save space the right child nodes are shown as down nodes. The three dots at the bottom indicate nodes, subsequences, to which the sieve is applied further. The tree dots on the top left indicate $EEE$ to which the sieve is applied.

After each application of the sieve to a parent's subsequence each member, number, is not missed and gets a new position in the left or right child's subsequences. Thus, after each finite number of iterative applications of the sieve every number is classified and occurs in one and only one \textit{terminal} subsequence. Since 1) the sieve is applied to every $E_m{E}$ subsequence creating its right child node $OE_m{E}$ from all numbers on odd positions, 2) the position of any number, except zero, alternates between even and odd after each new placement, and 3) the left nodes $EO_k{O}$ and $EOE_m{E}$ are created from each odd $O_k{O}$ and even $OE_m{E}$ subsequence, any number, except zero, after a finite number of iterative applications of the sieve becomes a member of a fixed subsequence. The numbers within a fixed subsequence are equivalent under the iterative sieve application.

Weather a subsequence is fixed or not its name is unique. It also associates with a unique formula of infinite arithmetic progression linearly expressing any number within the sequence from its position counted from zero. Given a number, a reverse linear function returns its position. The direct linear functions of $p$ are in the nodes of the tree one-to-one mapping the tree of names.
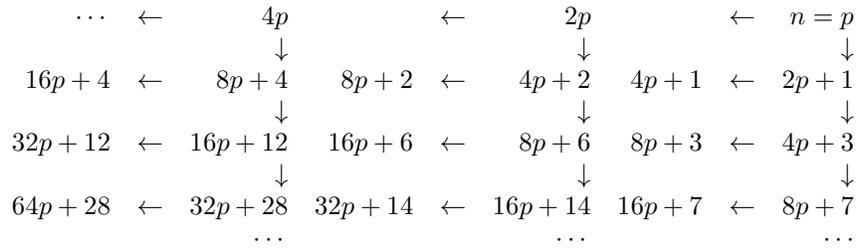
\begin{displaymath}
\begin{array}{rrrrrrrrr}
\cdots&\gets&4p&&\gets&2p&&\gets&n = p\\
&&\downarrow&&&\downarrow&&&\downarrow\\
16p+4&\gets&8p+4&8p+2&\gets&4p+2&4p+1&\gets&2p+1\\
&&\downarrow&&&\downarrow&&&\downarrow\\
32p+12&\gets&16p+12&16p+6&\gets&8p+6&8p+3&\gets&4p+3\\
&&\downarrow&&&\downarrow&&&\downarrow\\
64p+28&\gets&32p+28&32p+14&\gets&16p+14&16p+7&\gets&8p+7\\
&&\cdots&&&\cdots&&&\cdots
\end{array}
\end{displaymath}
Replacing the position $p$ in the formula of a parent's node with the even $2p$ or odd $2p+1$ position produces the formula for the left or right child nodes.

Geometrically, the points on the lines specified by the formulas associated with the fixed subsequences and computed for non-negative integer positions $p$ do not coincide. This is true for terminal sequences on each iteration.

Iterals of linear functions computed for $2p$ and $2p + 1$ initial values and different number of iterations together with their nesting properties handle complexity. In order to get the iteral expression for a node of the tree the following common procedure is available. It is based on already described iteral properties and the fact that in a tree a unique path connects any two nodes. Our paths begin from $\mathrm{N_0}$, which has the simplest linear function $p$. Making a move to a left node, apply \CYRI$_{2p}^{1}$ to a parent's function. Making a move to a right node, apply \CYRI$_{2p+1}^{1}$ to a parent's function. At the start up, as we have seen in the section on iteral expression, the two identities are obtained
\begin{displaymath}
E, 2p = \textrm{\LARGE{\CYRI}}_{2p}^{1}p \gets \mathrm{N_0} \to \textrm{\LARGE{\CYRI}}_{2p+1}^{1}p = 2p+1, O
\end{displaymath}
The iterals are nested on the path. As an example, consider the iterals nesting for the move $\mathrm{N_0}\to O \to OO \to EOO$
\begin{displaymath}
\textrm{\LARGE{\CYRI}}_{2p}^{1}\textrm{\LARGE{\CYRI}}_{2p+1}^{1}\textrm{\LARGE{\CYRI}}_{2p+1}^{1}p =
\textrm{\LARGE{\CYRI}}_{2p}^{1}\textrm{\LARGE{\CYRI}}_{2p+1}^{1}(2p+1) = 
\textrm{\LARGE{\CYRI}}_{2p}^{1}(4p+3) = 8p + 3.
\end{displaymath}
We shall be most interested in the fixed subsequences. They are subdivided on the subsequences $EO_k{O}$ containing only odd numbers
\begin{displaymath}
\textrm{\LARGE{\CYRI}}_{2p}^{1}
\textrm{\LARGE{\CYRI}}_{2p+1}^{k}(2p+1) = 2^{k+2}p+2^{k+1}-1, \; k = 0, 1, \; \ldots \; ,
\end{displaymath}
and subsequences $EO_l{E}_m{E}$
\begin{displaymath}
\textrm{\LARGE{\CYRI}}_{2p}^{1}\left(\textrm{\LARGE{\CYRI}}_{2p+1}^{1}\cdots
\textrm{\LARGE{\CYRI}}_{2p+1}^{1}\right)_{l}
\textrm{\LARGE{\CYRI}}_{2p}^{m}2p = 2^{m + l + 2}p + 2^{m+1}(2^l -1),
\end{displaymath}
\begin{displaymath}
m = 0, 1, 2, \; \ldots; \; l = 1, 2, 3, \; \ldots; \; p = 0, 1, 2, \; \dots \; .
\end{displaymath}
containing only even numbers. Zero always remains in a subsequence of even numbers with the growing name $EE_{\infty}E, \; l = 0, \; 2^{\infty+2}p, p=0$. By construction of the tree, we proved that each number becomes a member of one and only one of these subsequences. Thus, these iterals define the one-to-one correspondence between positive odd integers $a$ and pairs of nonnegative integers $(k, p)$ and between positive even integers $b$ and triplets of integers $(m \geq 0, l \geq 1, p \geq 0)$. In other words, $\forall$ odd $a \in \mathrm{N} \; \exists$ unique $k,p \in \mathrm{N_0}$ solving the equation
\begin{displaymath}
a = 2^{k+2}p+2^{k+1}-1.
\end{displaymath}
Similarly, $\forall$ even $b \in \mathrm{N} \; \exists$ unique $m,p \in \mathrm{N_0} \wedge l \in \mathrm{N}$ solving the equation
\begin{displaymath}
b =  2^{m + l + 2}p + 2^{m+1}(2^l -1).
\end{displaymath}
While $b = 0, \; l = 0, \; p = 0$ resolves the last equation for any $m$, we keep zero in $EE_{\infty}E$. Given a positive number $n$ the following algorithm consisting of three if-branches computes its $k$ and $p$ or $m, \; l$ and $p$.

If $n = 0$, then $l = 0, \; p = 0,$ and $m$ is arbitrary, set $m = \infty$.

If $n$ is odd, then subtract one and divide the result by two. If the result is odd, then repeat subtraction of one and division by two until an even number or zero is obtained. This is the stopping number. Count the number of repetitions on the way down to the stopping number: $k = counter_k - 1$. Divide the stopping number by two. The result is $p$. Example, $n = 39 \to 19 \to 9 \to 4$. We stop at an even number, or zero, after making three subtractions of one followed by division by two. Then, $k = 3 - 1 = 2$. The $p = \frac{stopping \; number}{2} = \frac{4}{2} = 2$. Verifying:
$2^{k+2}p+2^{k+1}-1 = 2^{2+2}2+2^{2+1}-1 = 39$. Example: $n = 7 \to 3 \to 1 \to 0; \; k = 3 - 1 = 2; \; p = \frac{0}{2} = 0$. Verifying: $2^{k+2}p+2^{k+1}-1 = 2^{2+2}0+2^{2+1}-1 = 7$. The numbers 7 and 39 are equivalent. They belong to one fixed subsequence $EOOO$ given by $16p + 7$. The number 5 corresponds to $(k = 0, p = 1)$ and belongs to $EO$ given by $4p + 1$. The numbers 5 and 7 are not equivalent since they belong to two different fixed classes. Of course, both are members of $O$ given by $2p + 1$. However, $O$ is not a fixed subsequence and a subject for application of the sieve.

If $n$ is even, then it can be divided $m + 1$ times by two giving the odd number $2^{l + 1}p + (2^l -1)$. The $m = counter_m - 1$. The odd number is unit, indicating that given $n$ is the power of two $= 2^{counter_m} = 2^{m+1}$, only, if $p = 0$ and $l = 1$. This means that all powers of two are extracted each by one and only one fixed subsequence of even numbers $EOE_m{E}$ and become their initial elements with positions $p = 0$. Thus, all powers of two are pair wise not equivalent. They are distributed between the fixed subsequences of even numbers as
\begin{displaymath}
 2 = 2^1 \in EOE : 8p + 2, \; 4 = 2^2 \in EOEE : 16p + 4, \; \ldots \; ,
\end{displaymath}
\begin{displaymath}
 2^{m + 1} \in EOE_m{E}:2^{m + 3}p+2^{m+1}, \; \ldots \;.
\end{displaymath}
For $n = 2^{m+1}$ the unique triplet is $(m = counter_m - 1, l = 1, p = 0)$. When the odd number is not the unit, then we treat it as in the second if-branch above for the odd numbers and compute the $counter_k$. The $l = counter_k$ and the $p = \frac{stopping \; number}{2}$. Remember that $k = counter_k - 1$. Example, $n = 28 \to 14 \to 7; \; m = 2 - 1 = 1; \; 7 \to 3 \to 1 \to 0; \; l = 3; \; p = \frac{0}{2} = 0$. Verifying: $2^{m + l + 2}p + 2^{m+1}(2^l -1) = 2^{1 + 3 + 2}p + 2^{1+1}(2^3 -1) = 64p + 28 = 28$. 28 is the initial member of the even fixed subsequence $EOOOEE$. The latter may not contain a power of two. Example, $n = 38 \to 19; \; m = 1 - 1 = 0; 19 \to 9 \to 4; l = 2; p = \frac{4}{2} = 2$.Verifying: $2^{m + l + 2}p + 2^{m+1}(2^l -1) = 2^{0 + 2 + 2}p + 2^{0+1}(2^2 -1) = 16p + 6 = 38$. Thus, $28 \in EOOE$. The latter may not have a power of two. Its initial member is $6$. 28 and 38 are not equivalent.

The unique mapping $0 \leftrightarrow (m = 0, l = 0, p = \infty), \; odd \leftrightarrow (k, p),$ and $even \leftrightarrow (m, l, p)$ and subdivision on equivalence fixed subsequences allows reordering $\mathrm{N_0}$ 
\begin{displaymath}
\cdots \prec E{O_l}E_m{E} \prec \cdots \prec EOE_m{E} \prec
\end{displaymath}
\begin{displaymath}
\cdots \prec E{O_l}EE \prec \cdots \prec EOEE \prec
\end{displaymath}
\begin{displaymath}
\cdots \prec EO_l{E} \prec \cdots \prec EOE \prec 0 \prec EO \prec \cdots \prec E{O_k}O \prec \cdots
\end{displaymath}
or with iterals
\begin{displaymath}
\cdots \prec
\textrm{\LARGE{\CYRI}}_{2p}^{1}\left(\textrm{\LARGE{\CYRI}}_{2p+1}^{1}\right)_{l}2^{m+1}p
\prec \cdots \prec
\textrm{\LARGE{\CYRI}}_{2p}^{1}\textrm{\LARGE{\CYRI}}_{2p+1}^{1}2^{m+1}p
 \prec
\end{displaymath}
\begin{displaymath}
\cdots \prec
\textrm{\LARGE{\CYRI}}_{2p}^{1}\left(\textrm{\LARGE{\CYRI}}_{2p+1}^{1}\right)_{l}2^2p
\prec \cdots \prec
\textrm{\LARGE{\CYRI}}_{2p}^{1}\textrm{\LARGE{\CYRI}}_{2p+1}^{1}2^2p
 \prec
\end{displaymath}
\begin{displaymath}
\cdots \prec
\textrm{\LARGE{\CYRI}}_{2p}^{1}\left(\textrm{\LARGE{\CYRI}}_{2p+1}^{1}\right)_{l}2p
\prec \cdots \prec
\textrm{\LARGE{\CYRI}}_{2p}^{1}\textrm{\LARGE{\CYRI}}_{2p+1}^{1}2p
 \prec
\end{displaymath}
\begin{displaymath}
0 \prec
\textrm{\LARGE{\CYRI}}_{2p}^{1}(2p+1)
\prec \cdots \prec
\textrm{\LARGE{\CYRI}}_{2p}^{1}
\textrm{\LARGE{\CYRI}}_{2p+1}^{k}(2p+1)
\prec \cdots
\end{displaymath}
All formulas for the iterals above are already derived.

The formulas for odd $a$, and even $b$ can be combined into one formula valid for any $n \in \mathrm{N}$, where $m = -1, 0, 1, \ldots, \; l = k + 1 = 1, 2, \ldots,$ and $p = 0, 1, \ldots$
\begin{displaymath}
n = 2^{m + l + 2}p + 2^{m + 1}(2^l - 1), \; n \leftrightarrow (m, l, p), \; m \in \mathrm{Z} \wedge m \geq -1, l \in \mathrm{N}, p \in \mathrm{N_0},
\end{displaymath}
\begin{displaymath}
n = 2^{m + 1}(2^{l+1}p + 2^l - 1).
\end{displaymath}
All natural odd numbers are produced for $m = -1$, all natural even numbers are generated with $m \in \mathrm{N_0}$, and all powers of two, a subset of even numbers, are obtained for $m \in \mathrm{N_0} \wedge l = 1 \wedge p = 0$. The formula does not miss any natural number and assigns a unique triplet $(m, l, p)$ for the latter. We also can get zero for $l = 0 \wedge p = 0$ and arbitrary $m$. Uniqueness does not hold for zero.

Waclaw Sierpinski calls Mersenne numbers the numbers represented as $2^l - 1, l = 1, 2, \; \ldots$ \cite[p. 76-79]{sierpinski}. The subset of prime numbers, Mersenne primes, are named after the Minimite friar, Farther, Marin Mersenne \cite[p. 334, 385]{boyer}, \cite[p. 11, 22]{goldman}. Computer specialists working with the standard binary representations of integers well know that $2^l - 1 = 1 \times 2^{l-1} + \; \ldots \; 1 \times 2^2 + 1 \times 2^1 + 1 \times 2^0 = \sum_{j = 0}^{j = l -1}2^j$. For instance, decimal $7$ is equal to binary $111 = 1 \times 2^2 + 1 \times 2^1 + 1 \times 2^0$. We see that each fixed arithmetic progression for odd numbers $a = 2^{k+2}p+2^{k+1}-1 = 2^{l+1}p + 2^l - 1 = 2^{l+1}p + M_l$ starts from a Mersenne number $M_l$ and has the progression difference equal to the power of two obtained after incrementing by one the position, starting from one, of a number within the sequence of Mersenne numbers. The formulas for the fixed arithmetic progressions of even numbers are obtained from the odd counterparts after multiplying by a power of two.

All spectrum of interesting mathematical questions related to the infinite subsequences can be applied to the equivalence arithmetic progressions considered in this section. For some of them, such as $EO : 4p + 1$, the representation of the prime numbers of the sequence by the sum of two squares have become classics of the number theory \cite[p. 51-53]{sierpinski}, \cite[p. 20]{goldman}. How the prime and/or Mersenne numbers are distributed between the fixed subsequences? It is worth to quote Sierpenski's words (author's translation to English) \cite[p. 80]{sierpinski}:

\textit{Question, whether a given infinite sequence, defined even in a simple way, contains infinitely many prime numbers, in general, is quite difficult.}

The task to show that iterals, iteral notation, and the sieve are natural means to derive the equivalence arithmetic progression formulas is completed.

The subdivision of $\mathrm{N}_0$ is interesting for studying the division by two applied to the even numbers in the even fixed subsequences on the left side of zero in conjunction with the operation $\frac{3x + 1}{2}$ applied to the odd numbers in the fixed subsequences in the right side of zero. This is done in the next sections.

\textbf{The operation $\frac{3x + 1}{2}$ applied to $EO_k{O}$.} There is unsolved convergence problem named after Lothar Collatz \cite{lagarias}. If the unit is given, then stop. If an even number greater than one is given, then divide it repeatedly by two until the unit or a greater odd number is obtained. If it is one, then stop. If the result is an odd number greater than one or an odd number is given initially, then apply to it repeatedly $\frac{3x + 1}{2}$. If an even number is obtained at any step, then switch to division by two. The unproved proposition is that for any positive  integer after the finite number of steps the result will be the unit. For instance, $3 \to 5 \to 8 \to 4 \to 2 \to 1$.

Within the scope of this article it is natural to denote the result of $\frac{3x + 1}{2}$ applied $k = 0, 1, \ldots$ times by the iteral \CYRI$_{n}^{k}\frac{3x + 1}{2}$. In order to get integer, the operation must be applied to an odd number $n = 2p + 1$. Under this condition, one iteration is equivalent to adding to an odd number a position of a next odd number within $O$. Indeed,
\begin{displaymath}
\textrm{\LARGE{\CYRI}}_{n = 2p + 1}^{1}\frac{3x + 1}{2} = \frac{6p + 4}{2} = 3p + 2 = 2p + 1 + p + 1 = n + p + 1.
\end{displaymath}
The result is an odd number, if $p$ is odd. Otherwise, the result is an even integer. We just proved that in terms of our subsequences, the operation converts an even-odd number to an even number: $EO \to E$. Getting an even number from an odd one is important for the entire proposition because this is the only way eventually to face with an even number, which is a power of two leading to final one. Unfortunately, the same operation converts an odd-odd number to an odd number. Let us prove that being applied $k + 1$ times to an odd number from the fixed subsequence of the order $k \; EO_k{O}$ an even number is obtained. Indeed, for the operation applied $j$ times to such a number we get
\begin{displaymath}
\textrm{\LARGE{\CYRI}}_{      \textrm{\LARGE{\CYRI}}_{2p}^{1}\textrm{\LARGE{\CYRI}}_{2p+1}^{k}(2p+1)        }^{j}\frac{3x + 1}{2} = \textrm{\LARGE{\CYRI}}_{2^{k + 2}p + 2^{k + 1} - 1}^{j}\frac{3x + 1}{2} = 
\end{displaymath}
\begin{displaymath}
 = 3^{j}2^{k + 2 - j}p + 3^{j}2^{k + 1 - j} - 1.
\end{displaymath}
As long as $j \leq k$ the result remains an odd number. When $j$ achieves $k + 1$ the result is even and equal to $3^{k+1}2p + 3^{k+1} - 1 = 3^{k+1}(2p + 1) - 1$. Thus, any odd number $n \leftrightarrow (k,p)$ after $k + 1$ iterations applying the operation $\frac{3x + 1}{2}$ is converted to an even number larger than $n$. It is larger because on each iteration it is added by increasing position of a next odd number.

What are $(k_{new}, p_{new})$ after application of $\frac{3x+1}{2}$ to an odd $n \leftrightarrow (k,p)$?
\begin{displaymath}
\textrm{\LARGE{\CYRI}}_{2^{k + 2}p + 2^{k + 1} - 1}^{1}\frac{3x + 1}{2} = 3^{1}2^{k + 1}p + 3^{1}2^k - 1 =
\end{displaymath}
\begin{displaymath}
 = (2 + 1)2^{k + 1}p + (2+1)2^k - 1 = 2^{k + 2}p + 2^{k + 1}p + 2^{k+1} + 2^k - 1.
\end{displaymath}
Let us apply the algorithm for odd numbers to find $k_{new}$ and $p_{new}$. The right most $-1$ is not a problem because we subtract one and then divide by two until even or zero stopping number. This will keep $-1$ at the end of expression intact. Thus, after $counter = k$ steps the stopping number is even $2^{2}p + 2p + 2 + 1 - 1 = 2^{2}p + 2p + 2$. Then, $k_{new} = k - 1$ and $p_{new} = 3p + 1$. This operation reduces the order of the odd fixed subsequence $EO_kO \to EO_{k-1}O$ and increases position within the new subsequence: $(k, p) \to (k-1, 3p+1)$. After $k$ iterations the converted number is still odd and belongs to $EO:4p+1$ with the order $0$. It is one step from becoming even. On $j$th application of the operation on  the way to smaller orders $\forall j \leq k$ the odd number is given by the pair $($\CYRI$_{k}^{j}(k-1),$ \CYRI$_{p}^{j}(3p+1)) = (k - j, 3^j{p}+\sum_{i=0}^{j-1}3^i) = (k - j, 3^j{p}+\frac{3^j-1}{2})$. The sum for $j = 0$ is set to zero as it is described in the historical section for $\sum$. Here are the starting and ending $(k_{last}, p_{last})$ nodes on the iteration way, when the numbers are still odd
\begin{displaymath}
(k, p) \in EO_k{O}:2^{k + 2}p + 2^{k + 1} - 1 \to (0, 3^{k}p+\frac{3^k-1}{2}) \in EO:4p+1.
\end{displaymath}
Single application of the operation strips off one $O$ from the name of an odd fixed subsequence. The $(k+1)$th application returns an even number $3^{k+1}(2p + 1) - 1$. It can be obtained as we did it above or from $4p+1$ replacing $p$ by $p_{last}$ and applying the operation the last $(k+1)$th time, when it still can be applied.

\textbf{The operation $\frac{x}{2}$ applied to $EO_l{E_m}E$.} The even numbers $(m, l, p)$ from these subsequences are given by $2^{m + l + 2}p + 2^{m+1}(2^l -1)$ and can be divided $(m+1)$ times by two
\begin{displaymath}
\textrm{\LARGE{\CYRI}}_{2^{m + l + 2}p + 2^{m+1}(2^l -1)}^{m+1}\frac{x}{2} = 2^{l + 1}p + 2^l -1.
\end{displaymath}
The result is an odd number $(k, p)$ from a subsequence $EO_k{O}$. Indeed, for $l = k + 1$ we get $2^{k + 2}p + 2^{k+1} -1$. Being applied $(m+1)$ times to a number from an even fixed subsequence, division by two strips off all suffixes $E$: $EO_l{E_m}E \to EO_l = EO_{l-1}O = EO_k{O}$. The original even number and the obtained odd number have the same numerical position each within its own subsequence. The intrinsic order of an odd number "hidden" in the structure of an even number does not change during this transformation. With respect to division by two $l = k + 1$ and $p$ are invariant. 

Using the current layout of the tree diagram, where the right child nodes are shown as the bottom nodes, division by two starting from an even fixed subsequence is a horizontal move to the right, sequentially visiting nodes of other fixed even subsequences with incrementally reduced number of $E$ in the suffix of the name. This continues until the node of an odd fixed subsequence with the same number of $O$s is visited. Then control is given to $\frac{3x+1}{2}$, which in $k$ vertical moves visits $EO$. The latter is a subject of the last $(k+1)$th application giving an even $3^{k+1}2p + 3^{k+1} - 1$. Which even subsequence $EO_{l_{new}}{E_{m_{new}}}E : 2^{m_{new} + l_{new} + 2}p_{new} + 2^{m_{new}+1}(2^{l_{new}} -1)$ does it belong to?
\begin{displaymath}
\begin{array}{rcrcc}
&&&EO&\\
&&\swarrow 1,\; \frac{3x+1}{2}&&\\
&EO_?E_?E = EO_{l_{new}}{E_{m_{new}}}E&&\uparrow&k = l - 1,\\
&&&&\frac{3x+1}{2}\\
E{O_k}O{E_m}E&\to&&E{O_k}O&\\
&m+1, \; divisions \; by \; 2&&&
\end{array}
\end{displaymath}

We have proved earlier that $2^{m+1} \in EO{E_m}E : 2^{m + 3}p+2^{m+1}$ is the initial $p = 0$ member of the subsequence. In order the main conjecture would be true, eventually, a move from $EO$ should visit one of these subsequences and exactly its initial member: $l_{new}=1 \wedge p_{new}=0$. "Eventually" emphasizes that this does not have to happen on the first loop. Then, both endless 1) cycling between $EO$ and non-initial member of $EO{E_m}E$ and 2) general loops must be excluded.

\textbf{The fate of even number $3^{k+1}2p + 3^{k+1} - 1$.} The even if-branch of the algorithm determining $(m, l, p)$ begins from division by two until an odd number is obtained. If the odd is one, then the source is the power of two, otherwise, $m = counter_{m} - 1$ and the odd if-branch of the algorithm takes control. The $k$ and $p$ are arbitrary and independent because any number can be given at the beginning. Without loosing generality, there are four possibilities: 1) $k = 0 \wedge p = 0$, then the even number is two and we do not have to be here because the original given number is $1 \in EO:4p+1$, indeed, $2 = \frac{3 \times 1 + 1}{2}$; 2) $k = 0 \wedge p = odd_p2^q, \; odd_p \geq 1, \; q \geq 0$, then the number is $6p + 2 = 3odd_p2^{q+1}+2$ and $\forall q>0 \; m_{new} = 0$; 3) $k = odd_k2^r \wedge p = 0, \; odd_k \geq 1, \; r \geq 0$, then the number is $3^{k+1} - 1$; and 4) $k = odd_k2^r \wedge p = odd_p2^q$ producing the most complex expression, where we apply the Newton's binom \cite{korn} for the second occurrence of $3^{k+1} = (2 + 1)^{k + 1}$. Euler's $\sum$ and Kramp's $!$ are in use
\begin{displaymath}
3^{k+1}2p + (2+1)^{k+1} - 1 =3^{k+1}odd_p 2^{q+1} + \sum_{j=0}^{k+1}\frac{(k+1)!}{j!(k+1-j)!}2^{k+1-j} - 1 =
\end{displaymath}
\begin{displaymath}
3^{k+1}odd_p 2^{q+1} + \sum_{j=0}^{k}\frac{(k+1)!}{j!(k+1-j)!}2^{k+1-j} =
\end{displaymath}
\begin{displaymath}
3^{k+1}odd_p 2^{q+1} + \sum_{j=0}^{odd_k 2^r}\frac{(odd_k 2^r+1)!}{j!(odd_k 2^r+1-j)!}2^{odd_k 2^r+1-j}.
\end{displaymath}
The $(k+1)$th term $1$ is taken out of the $\sum$ and the last value of the index $j$ is reduced to $k$. The $(k)$th term of the sum is equal to $(k + 1)2 = (odd_k 2^r + 1)2^1$. Comparing with other terms, for $r>0$ it contains the smallest factor represented by a positive power of two. For $r>0$ the sum is divisible by two only once. By the same reasoning in the third case for $r>0$ the number $3^{k+1}-1$ is divisible by $2$ only one time. Thus, $\forall odd_p, odd_k, q,r \in \mathrm{N} \; m_{new}=0$ and $3^{k+1}2p + 3^{k+1} - 1 \in EO_{l_{new}}E=EO_{k_{new}+1}E:2^{l_{new}+2}p_{new}+2(2^{l_{new}}-1)=2^{k_{new}+3}p_{new}+2(2^{k_{new}+1}-1)$.
After single division by two we get for the second, third, and fourth cases
\begin{displaymath}
3odd_p2^{q}+1, \; \textrm{odd for} \; q > 0 \; \textrm{and even for} \; q = 0;
\end{displaymath}
\begin{displaymath}
\sum_{j=0}^{odd_k 2^r}\frac{(odd_k 2^r+1)!}{j!(odd_k 2^r+1-j)!}2^{odd_k 2^r-j}, \; \textrm{odd for} \; r > 0 \; \textrm{and even for} \; r = 0;
\end{displaymath}
\begin{displaymath}
3^{k+1}odd_p 2^{q} + \sum_{j=0}^{odd_k 2^r}\frac{(odd_k 2^r+1)!}{j!(odd_k 2^r+1-j)!}2^{odd_k 2^r-j}.
\end{displaymath}
In the fourth case the number is odd for $(q>0 \wedge r>0) \vee (q=0 \wedge r=0)$ and even for
$(q>0 \wedge r=0) \vee (q=0 \wedge r>0)$. Determination of $m_{new}$ is not over for even result. While for any given natural number $n$ the algorithm returns $(m \geq 0, l \geq 1, p \geq 0)$ in a straight forward computation, without an explicit expression for $EO_{l_{new}}{E_{m_{new}}}E : 2^{m_{new} + l_{new} + 2}p_{new} + 2^{m_{new}+1}(2^{l_{new}} -1)$ it is difficult to move forward. The conjecture remains unproved. It is proved that any natural number, which is not a power of two, will be converted to an even number using the finite number of divisions by two and $\frac{3x+1}{2}$ operations, which can be determined in advance by setting the correspondence between the given number and pair $(k, p)$ or triplet $(m, l, p)$. Application of iterals to these operations looks natural.

The author is grateful to Florian Cajori, who many years ago has written the words suitable for ending this article \cite[p. 77]{cajori}

\textit{We know that intellectual food is sometimes more easily digested, if not taken in the most condensed form. It will be asked, To what extent can specialized notations be adopted with profit? To this question we reply, only experience can tell. It is one of the functions of the history of mathematics to record such experiences.}

\paragraph{Acknowledgments.}  I would like to thank Alexander Tumanov and a referee for reviewing the article and useful advices and criticism. Due to them I have added the second version of the iteral notation, where in ambiguous cases a variable of a function is specified explicitly.

\section{Appendix}
The author has written a program \textit{oneness}. The C++ code can be requested by email. Given a positive integer and radix of the system (decimal, octal, binary, hexadecimal, and, in general, from two to 36) it returns intermediate numbers for a Collatz's task. The program works with the numbers not exceeding in computations the Mersenne number $2^{31} - 1 = 2147483647$. Here is the output for the run: oneness 9 3
\begin{verbatim}
 0  9 IN  100 EO : 4p + 1 : (k=0, p=2) : 9
 1 14 3X  112 EO3E : 32p + 14 : (m=0, l=3, p=0) : 14
 2  7 D2   21 EO2O : 16p + 7 : (k=2, p=0) : 7
 3 11 3X  102 EO1O : 8p + 3 : (k=1, p=1) : 11
 4 17 3X  122 EO : 4p + 1 : (k=0, p=4) : 17
 5 26 3X  222 EOE : 8p + 2 : (m=0, l=1, p=3) : 26
 6 13 D2  111 EO : 4p + 1 : (k=0, p=3) : 13
 7 20 3X  202 EOE1E : 16p + 4 : (m=1, l=1, p=1) : 20
 8 10 D2  101 EOE : 8p + 2 : (m=0, l=1, p=1) : 10
 9  5 D2   12 EO : 4p + 1 : (k=0, p=1) : 5
10  8 3X   22 EOE2E : 32p + 8 : (m=2, l=1, p=0) : 8
11  4 D2   11 EOE1E : 16p + 4 : (m=1, l=1, p=0) : 4
12  2 D2    2 EOE : 8p + 2 : (m=0, l=1, p=0) : 2
13  1 D2    1 EO : 4p + 1 : (k=0, p=0) : 1
\end{verbatim}
The first column is the step counted from zero. The second column is an intermediate number, where nine on the line zero is the input. The third column is the name of the operation, where IN is input, and 3X or D2 is $\frac{3x+1}{2}$ or division by two applied to a previous number. The fourth column is the number represented with the radix three requested in the command line. The last descriptor with the colon delimiters contains the name of the subsequence, the formula of the arithmetic progression, the pair $(k, p)$ or triplet $(m, l, p)$, and repetition of the same number computed after reusing the pair or triplet to make sure that the program works correctly.

\bigskip

\noindent\textbf{Valerii Salov} received his M.S. from the Moscow State University, Department of Chemistry in 1982 and his Ph.D. from the Academy of Sciences of the USSR, Vernadski Institute of Geochemistry and Analytical Chemistry in 1987.  He is the author of the articles on analytical, computational, and physical chemistry, the book Modeling Maximum Trading Profits with C++, \textit{John Wiley and Sons, Inc., Hoboken, New Jersey}, 2007, and papers in \textit{Futures Magazine}.

\noindent\textit{v7f5a7@comcast.net}

\end{document}